\documentclass[12pt]{article}
\usepackage{eepic,amsmath,amsthm,amssymb} 
\usepackage[all]{xy}

\newcommand{\ds}{\displaystyle} 

\newcommand{\db}{\bar{\partial}}
\newcommand{\dc}{\check{\delta}}

\newcommand{\im}{\operatorname{Image}}
\renewcommand{\ker}{\operatorname{Ker}}
\newcommand{\der}{\operatorname{Der}}
\newcommand{\sym}{\operatorname{\mathcal{S}}}

\renewcommand{\L}{\mathfrak{L}} 
\newcommand{\g}{\mathfrak{g}}
\newcommand{\h}{\mathfrak{h}} 
\newcommand{\w}{\mathrm{wt}}
\newcommand{\A}{\mathfrak{A}}
\newcommand{\F}{\mathfrak{F}}
\newcommand{\C}{\mathcal{C}} 
\newcommand{\R}{\tilde{R}}
\renewcommand{\P}{\mathbf{P}} 
\newcommand{\M}{\mathcal{M}}
\newcommand{\N}{\mathcal{N}} 
\renewcommand{\O}{\mathcal{O}}
\renewcommand{\H}{\mathcal{H}}

\theoremstyle{plain} 
\newtheorem{theorem}{Theorem}[section]
\newtheorem*{theorem*}{Theorem} 
\newtheorem{lemma}[theorem]{Lemma} 

\newtheorem{corollary}[theorem]{Corollary}
\newtheorem*{d-d}{$\mathbf{\Delta-\db}$ Lemma}
\newtheorem*{corollary*}{Corollary} 

\theoremstyle{definition}
\newtheorem{definition}{Definition} 
\newtheorem*{remark*}{Remark}
\newtheorem{remark}{Remark}

\title{Hypersurfaces and generalized deformations} 
\author{John Terilla}
\date{\today}

\begin{document}

\maketitle

\abstract{The moduli space of generalized deformations of a
  Calabi-Yau hypersurface is computed in terms of the Jacobian ring of
  the defining polynomial. The fibers of the tangent bundle to this
  moduli space carry algebra structures, which are identified using
  subalgebras of a deformed Jacobian ring.

\section{Introduction}

\subsection{Background}

In his landmark paper \cite{W3} ``Mirror Manifolds and Topological
Field Theories,'' Witten discussed his $B(X)$ model topological field
theory (TFT) associated to a Calabi-Yau manifold $X$ and introduced
the problem of deforming $B(X)$ as a TFT.  Denote by $\M(X)$, or just
$\M$, the moduli space of deformations of $B(X)$. Witten concluded
that the deformations of $B(X)$ are parametrized by
$H^{\bullet}(\Theta_X^\bullet)$.  Not surprisingly, the deformations
of the complex structure of $X$---these correspond to elements of
$H^{1}(\Theta_X)\subseteq H^{\bullet}(\Theta_X^\bullet)$---are among
the deformations of $B(X)$.  For this reason, elements of $\M$, are
called generalized, or extended, deformations of $X$.

It is natural to look for a mathematical deformation problem (inspired
by, but without using, physical field theory) that gives rise to $\M$
as a moduli space.  One is led into formal deformation theory \cite{S}
and its Lie algebraic expression \cite{SS2,SS,GM}.  Following Witten,
Ran \cite{R} wrote about an algebraic deformation problem associated to
a Calabi-Yau manifold $X$ whose first order deformations are in one to
one correspondence with $H^{\bullet}(\Theta_X)$.

Independently, Boris Dubrovin \cite{dub} studied deformations of
abstract topological field theories and discovered that their moduli
spaces always carry additional structures, which he axiomatically
treated and named \emph{Frobenius manifolds}. Barannikov and
Kontsevich \cite{BK} considered a deformation problem, different from
Ran's, with moduli space $H^{\bullet}(\Theta_X^\bullet)$ and proved
that this space is a Frobenius manifold.  The Barannikov and
Kontsevich theory is connected via the formality theorem \cite{K1} to
the moduli of  $A_\infty$ deformations of the category of coherent
sheaves on a complex manifold $X$ introduced in \cite{K2}.

\subsection{Contents}

In the case that $X$ is a Calabi-Yau projective hypersurface, one
expects that the extended deformation theory, the moduli space $\M$,
and some of the additional structures on $\M$ will have algebraic
descriptions in terms of the defining polynomial $f$.  We give these
algebraic descriptions in this paper.

In section \ref{background}, we review the Barannikov-Kontsevich
construction of $\M(X)$ for an arbitrary Calabi-Yau $X$ and recall
some of the structures implied by $\M$'s status as a Frobenius
manifold.  In particular, each fiber of $T\M$ has an associative
product.

In section \ref{answer}, we discuss the case in which $X$ is a
hypersurface with defining polynomial $f$. The main results of the
paper are stated in this section.  We describe $\M(X)$ in terms of a
subvectorspace of the Jacobian algebra of $f$. Then, a (primitive)
element $a\in \M$ can be identified with a polynomial $g$ and the
algebra structure on $T_a\M$ is identified with a subalgebra of the
Jacobian ring of $f+g$. There are similarities here to the deformation
theory of a function and to the complex deformations of a projective
hypersurface.  However, we know of no other place where the specific
algebras described here have appeared.
  
The remaining sections are devoted to constructing a deformation
theory, equivalent to the one introduced in \cite{BK}, that yields the
polynomial descriptions presented in section \ref{answer}.
  
The work in \cite{BK} relies on a Dolbeault resolution of
$H^\bullet(\Theta^\bullet_X)$.  Section \ref{cech} lays the foundation
for using \v Cech cochains instead.  Even though section \ref{cech} is
used in this paper only as a means to verify the presentation in
section \ref{answer}, we consider the successful replacement of the
Dolbeault complex by the \v Cech complex to be a main result
of the paper.

Section \ref{L} describes an algebraic model for $\Theta_X$ in the
case that $X$ is a hypersurface.  We describe a sheaf of differential
graded Lie algebras $(\L,d_f,[\,,\,])$ that is quasi-isomorphic to the
sheaf of holomorphic vector fields on $X$.  The achievement of this
model is that it describes $\Theta_X$ as the cohomology of a simple
complex where all the information about $X$ is carried in the
differential.  This makes the model well suited to being deformed.

In section \ref{F}, we use $\L$ to form $(\C,D)$---the complex of \v
Cech cochains with values in powers of $\L.$ This complex is our
analogue of the Dolbeault resolution of $H^\bullet(\Theta_X^\bullet)$
used in \cite{BK}.  We then give the complex $(\C,D)$ the structure of
an $L_\infty$ algebra.  The theorems stated in section \ref{answer}
follow from studying the deformation theory controlled by this
$L_\infty$ algebra.

Concluding remarks are contained in section \ref{conc}.

\section{A picture of $\M$}\label{background}

\subsection{$\M$ as a moduli space of generalized deformations}\label{M1}
Let $X$ be a Calabi-Yau manifold, let $\Theta_x$ denote the
holomorphic tangent sheaf of $X$, and let $\Theta_X^\bullet$ denote
$\ds \oplus_{p=0}^{\text{dim}(X)} \Lambda^p \Theta_X$.  Consider the
Dolbeault resolution $(\g,\db)$ of $H^\bullet(\Theta^\bullet_X):$
\begin{equation*}\label{g}
\ds \g=\bigoplus_{k}\g^k, \quad \g^k=\bigoplus_{k=p+q-1}\Gamma\left(\Theta^p_X\otimes
  \overline{\Omega}^q_X\right).
\end{equation*}
Throughout this paper, $\g$ will denote the graded vector space
defined above.  Note that $H^\bullet(\Theta^\bullet_X)=H(\g,\db)$, as
the cohomology of a graded vector space, inherits a grading.  With
this grading, a homogeneous element $\alpha \in H^q(\Theta^p_X)$ has
degree $p+q-1.$ We will not view
$H^\bullet\left(\Theta_X^\bullet\right)$ as a bigraded object.

There is a Schouten bracket $[\,,\,]$ on $\g$ extended from the
ordinary bracket of vector fields which makes $(\g,\db,[\,,\,])$ into
a differential graded (dg) Lie algebra.

It is a general fact that for any dg Lie algebra $(\h,d,[\,,\,])$,
there is a deformation theory controlled by $(\h,d,[\,,\,])$,
represented by a moduli space.  Informally, the differential $d$ is
deformed and the moduli space consists of (equivalence classes of)
elements $\alpha$ making
$$d_\alpha\overset{\mathrm{def}}{=}d+\mathrm{ad}(\alpha)$$
into a
differential.  The deformed differential $d_\alpha$ may contain
parameters (from the maximal ideal $m$ of an Artin local ring: $\alpha
\in \h \otimes m$).  The map $d_\alpha$ is always a derivation and the
condition that $d_\alpha^2=0$ translates into the condition on
$\alpha$:
$$d\alpha+\frac{1}{2}[\alpha,\alpha]=0,$$
which is sometimes called
the Mauer-Cartan equation, or the master equation.

More formally, to any Artin local ring $A$ with maximal ideal $m$, one
defines a functor $Def_\h$ from Artin local algebras with residue field $k$
to sets by
$$Def_\h(A)=\{\alpha \in (\h\otimes m)^1:
d\alpha+\frac{1}{2}[\alpha,\alpha]=0\}/\sim$$
where $\sim$ is the
equivalence generated infinitesimally by $(\h\otimes m)^0$, which acts
on $\alpha \in (\h\otimes m)^1$ by $\beta \cdot \alpha = [\beta,\alpha]-d\beta$. 
One may extend the functor $Def_\h$ to one acting
on \emph{graded} local Artin rings, which we denote by
$Def^\mathbb{Z}_\h$.  In \cite{S}, the problem of representing a
functor of Artin rings is discussed. The functor $Def_\h$ (or
$Def^\mathbb{Z}_\h$), is said to be represented by $\O$, a projective
limit of local Artin rings, if $Def_\h(A)=Hom(\O,A)$ (or
$Def^\mathbb{Z}_\h(A)=Hom(\O,A)$).  In the case that $\O=\O_\N$ is the
ring of local functions of a pointed space $\N$, then $\N$ is called
the moduli space for $Def_\h$ (or $Def^\mathbb{Z}_\h$).

We refer the reader to \cite{SS2,SS,GM,BK,Ba,K1} for more details, or
to \cite{K1} or section \ref{deformation_theory} of this paper, where
a more general deformation theory governed by an $L_\infty$ algebra is
discussed.

\begin{definition}\label{M}
  Define $\M$ to be the moduli space for $Def^\mathbb{Z}_\g.$
\end{definition}

\subsection{Algebraic structures on $\M$}
We now review two theorems which appear in \cite{BK}.

The set $Def^\mathbb{Z}_\h(k[t]/t^2)$ is the set of infinitesimal deformations.  
It has a natural vector space structure which is isomorphic to $H(\h)$.  
Provided the moduli space $N$ 
exists, the infinitesimal deformations are naturally isomorphic to 
$\hom(\O_N, k[t]/t^2)$ which can be identified with  $T_0 N$.
However, the moduli space may not exist as a smooth manifold and
infinitesimal deformations may be obstructed.  But the dg Lie algebra 
$\g$ is special:

\begin{theorem}\label{formal}
  $\M$ can be identified with a formal neighborhood of zero in the
  graded vector space $H(\g,\db)$.
\end{theorem}

Now there is some additional algebraic structure on $\M$ which we describe.
The associative $\wedge$ product on
$\Gamma\left(\overline{\Omega}^q_X\right)$ brings an additional
multiplication making $(\g,\db,[\,,\,],\wedge)$ a differential
Gerstenhaber algebra.  This means that $(\g,\db,[\,,\,])$ is a
differential graded Lie algebra, that $(\g,\db,\wedge)$ is a
differential graded commutative associative algebra (after a degree
shift), and that these structures are compatible via the odd Poisson
identity:
\begin{equation*}\label{gerst}
[\alpha, \beta \wedge \gamma]=[\alpha,\beta]\wedge
\gamma + (-1)^{(i+1)j} \beta \wedge[\alpha,\gamma] \quad \text{for
  }\alpha \in \g^i, \beta \in \g^j, \gamma \in \g.
\end{equation*}

Let $a \in H(\g,\db).$ Viewed as a deformation of $\db$, there
corresponds a differential $\db_\alpha$ with $[\alpha]=a.$
The compatibilities among $[\,,\,],$ $\db,$ and $\wedge$ ensure that
the shifted cohomology $H(\g,\db_\alpha)$ inherits an associative
product from $(\g,\wedge)$.

\begin{theorem}\label{alg}
  For all $\alpha \in \M$, $H(\g,\db_\alpha) \simeq H(\g,\db)$ as
  vector spaces.
\end{theorem}

Note that unless the degree of $\alpha$ is $+1$, the differential
$\db_\alpha$ will not be a homogeneous derivation and hence
$H(\g,\db_\alpha)$ will not be graded.  The isomorphism in theorem
\ref{alg} is only as ungraded vector spaces.

Since, as vector spaces, $H(\g,\db_\alpha) \simeq H(\g,\db)$ and
$H(\g,\db_\alpha)$ is an associative algebra, theorems \ref{formal}
and \ref{alg} may be assembled to view $\M$ as a (formal neighborhood
of zero in a) vector space with an associative product on each fiber
of the tangent bundle $T\M$: simply identify $T_{a}\M$ with
$H(\g,\db_\alpha)$.

\subsection{$\g$ as a dGBV algebra}
An effective setting in which to explain these two theorems is to
regard $\g$ as a differential Gerstenhaber-Batalin-Vilkovisky (dGBV)
algebra.  See the remark at the end of section 7 in \cite{BK} and the
elaboration in \cite{M1}.

The BV suffix is added because there is another degree one
differential $\Delta:\g\rightarrow \g$ that commutes with $\db$, is a
derivation of the bracket, and, together with $\wedge$, ``generates''
the bracket:
\begin{equation}\label{BV}
(-1)^i \Delta(\alpha\wedge \beta)- (-1)^i \Delta(\alpha)\wedge \beta -
\alpha\wedge\Delta(\beta)= [\alpha,\beta]
\text{ for $\alpha \in \g^i,\,\beta \in\g.$}
\end{equation}

This behavior of $\Delta$ together with the fact that $(\g,\db,
[\,,\,],\wedge)$ is a Gerstenhaber algebra, makes the quintuple
$(\g,\db,\Delta,[\,,\,],\wedge)$ what is called a dGBV algebra (there
is a nice survey of dGBV algebras in \cite{VS}).  Furthermore, there
is an additional homological property relating the two differentials
that make the particular dGBV algebra $(\g,\db,\Delta,[\,,\,],\wedge)$
extraordinary:

\begin{d-d}\label{dd}
  $\mathrm{Im}(\db)\cap\mathrm{Ker}(\Delta)=
  \mathrm{Im}(\Delta)\cap\mathrm{Ker}(\db)=\mathrm{Im}(\db\Delta).$
\end{d-d}

\begin{corollary}\label{formal2}
  The dg Lie maps, $$(\ker(\Delta),\db,[\,,\,]) \rightarrow
  (\g,\db,[\;,\;]) \text{ and } (\ker(\Delta),\db,[\,,\,]) \rightarrow
  (H(\g,\Delta),0,0)$$
  induce isomorphisms in cohomology.
\end{corollary}

Any map of dg Lie algebras inducing an isomorphism in cohomology
provides an isomorphism of deformation functors
\cite{SS2,SS,GM,K1,BK,Man}.  Corollary \ref{formal2}, then, implies that
the moduli space associated to $(\g,[\,,\,],\db)$ is equal to the
moduli space associated to $(H(\g,\Delta),0,0)$.  Because both the
differential and the bracket in $H(\g,\Delta)$ vanish, the moduli
space associated to $(H(\g,\Delta),0,0)$ can be simply identified with
a formal neighborhood of zero in the graded vector space
$H(\g,\Delta)$.  So, $\M$, defined as the moduli space associated to
$(\g,[\,,\,],\db)$, can be identified with $H(\g,\Delta).$ 
The $\db-\Delta \text{ Lemma}$ reveals that
$H(\g,\Delta){\simeq} H(\g,\db)$ by establishing natural
isomorphisms between both cohomologies and the quotient
$$\left (\ker(\Delta)\cap \ker(\db)\right)/ \left(\im \Delta \db \right).$$
Given $a\in M$, there corresponds an $[\alpha]\in H(\g
,\Delta)$ with $\Delta(\alpha)=\db(\alpha)=0.$  The identification
of $a\in M$ with $[\alpha]\in H(\g,\db)$ gives Theorem \ref{formal}.

It can be checked that if $\Delta(\alpha)=0$ then $\db-\Delta \text{
  Lemma }\Rightarrow \db_\alpha-\Delta \text{ Lemma.}$  
So, Theorem
\ref{alg} follows from having established that each element in $a\in\M$
can be represented by an element $[\alpha]\in H(\g,\db)$ with $\alpha\in \ker(\Delta)$.  We have
quasi-isomorphisms
$$
\xymatrix{ & \ar@{<.}[dl]_{\db-\Delta \text{ lemma }\quad}
  (H(\g,\Delta),0,0)
  \ar@{<.}[dr]^{\quad \db_\alpha-\Delta \text{ lemma}} \\
  (\g,\db_,[\;,\;]) & & (\g,\db_\alpha,[\;,\;]) }
$$
Therefore, as (ungraded) vector spaces $H(\g,\db_\alpha)\simeq
H(\g,\db)$, since they are both isomorphic, as vector spaces, to
$H(\g,\Delta)$.  Note that this diagram does not imply a
stronger isomorphism between $H(\g,\db_\alpha)$ and $H(\g,\db)$;
$\Delta$ is not a derivation of $\wedge$ and hence $H(\g,\Delta)$ is
not an associative algebra---the arrows are dgLie maps.

\section{Hypersurfaces}\label{answer}
Let $S=\mathbb{C}[x^0,\ldots,x^{n}]$. If $g\in S$ is a homogeneous
polynomial, denote its polynomial degree by $\w(g)=k$ and call it the
\emph{weight} of $g.$ So, $S=\oplus_k S^{[k]}$, where $S^{[k]}=\{g\in
S: \w(g)=k\}$.  Suppose $f\in S^{[\nu]}$ is nonsingular.  Denote
$\frac{\partial f}{\partial x^i}$ by $f_i$ and let $J_f$ be the ideal
of $S$ generated by $\{f_i\}_{i=0}^{n}.$ We have a subalgebra of
$S/J_f$ defined by
$$R=\im \left( \oplus_k S^{[k\nu]}\rightarrow S/J_{f}\right).$$
Note
that $J_f$ is a homogeneous ideal in $S$ so $R$ (like $S$) is
weighted: $R=\oplus_k R^k$ where $R^k:=\im(S^{[k\nu]}\rightarrow
S/J_f)$.  Let $$X=\{x\in\P^n:f(x)=0\}.$$   From now on, we assume $X$ 
is a Calabi-Yau manifold.  This
assumption means that $n+1=\nu$ and, as we already assume, $f$ is nonsingular.  So, $X$ is a degree $n+1$ hypersurface of
dimension $n-1$ in $\P^n.$

For $g\in \oplus_{k}S^{[k\nu]}$, define the deformed algebra
$$R_{f+g}=\im\left( \oplus_k S^{[k\nu]}\rightarrow S/J_{f+g}\right)$$
where $J_{f+g}$ is the ideal generated by $\left\{\frac{\partial
    (f+g)}{\partial x_i}\right\}$.  Note that $f+g$ is not, in
general, homogeneous and so $R_{f+g}$ is not weighted.

\begin{definition}Define a graded $\mathbb{C}$ algebra $\R$ and an
  ungraded $\mathbb{C}$ algebra $\R_{f+g}$ by
\begin{gather*}
  \R=R\oplus \mathbb{C}e_0 \oplus \cdots \oplus \mathbb{C}e_{n-1}
  \intertext{and} \R_{f+g}=R_{f+g}\oplus \mathbb{C}e_0 \oplus \cdots
  \oplus \mathbb{C}e_{n-1}
\end{gather*} 
with products extended from the products in $R$ and $R_{f+g}$ by
$$e_i\cdot e_j=e_i \cdot [h]=0 \text{ for all }i,j=0,\ldots n-1 \text{
  and }h\in S^{[k\nu]}, k>0.$$
The grading on $\R$ is defined by $\deg([h])=2k$ for $h\in S^{[k\nu]}$ and $\deg([e_i])=n-1.$
\end{definition}

We now state the main theorems (proved in sections \ref{cech} and \ref{theta}) in this paper.  They follow from analyzing the deformation theory governed by an $L_\infty$ algebra $(\C,Q^\C)$ that we construct (cf. definition \ref{qi})).

\begin{theorem}\label{result1}
$(\C,Q^\C)$ is quasi-isomorphic to $(\g,\db,[\,,\,])$.
Hence, $\M(X)$ can be identified with a formal neighborhood of zero in the vector space 
$\R$.
\end{theorem}

Our $L_\infty$ algebra carries an associative product
compatible with $D$, the differential term of $Q^\C$, and 
the cohomology of $(\C,Q^\C)$ is  naturally isomorphic to $\R$ as an algebra.  
An analysis of the shifted  cohomology rings, yields:

\begin{theorem}\label{result2} For each $[\alpha] \in \M(X)$ there
  corresponds a differential $D_{[\alpha]}$, a deformation of $D$, 
  and a polynomial $g \in S$ so that
  as (ungraded) associative algebras $H(\C,D_{[_\alpha]})\simeq
  \R_{f+g}.$
\end{theorem}

We make a remark about each theorem.

\begin{remark}
  Since $$\oplus_{p,q}H^{q}(\Theta^p_X)\simeq
  \oplus_{p,q}H^{q}(\Omega_X^{n-1-q})\simeq
  \oplus_{p,q}H^{p,q}(X)\simeq \oplus_i H_{DR}^i(X),$$
  all of the
  tangent cohomology of $\oplus_{p,q}H^{q}(\Theta^p_X)$ appears as
  DeRham cohomology of $X$. Using Hodge theory and residues \cite{G1}
  one can identify the primitive (middle) cohomology
  $H^{n-1-k,k}(X)$ with $R^k.$ This accounts for the ``primitive'' tangent
  cohomology contained in $\oplus_k H^k(\Theta^k).$  The rest of the
  tangent cohomology is inherited from projective space as $\oplus_k
  H^{k}(\Theta_X^{n-k-1})$ and each inherited piece
  $H^{k}(\Theta_X^{n-1-k})$ is one dimensional, adjoined to $R$ as
  $\mathbb{C}e_k.$ It is not surprising that
  $\tilde{R}\simeq H^\bullet(\Theta^\bullet_X)\simeq H(\g,\db)$, and the
  reader may recognize the grading on $\R$ as the grading acquired
  from an identification of $\R$ with $H^\bullet(\Theta^\bullet_X)$.  
  However, in this
  paper,  $\R$ arises as the cohomology of a graded algebra and as a 
  moduli space---theorem \ref{result1} is
  revealed in a natural setting.
  \end{remark}

\begin{remark} Theorem \ref{result2}, is stated in terms of the fairly simple algebras $\R_{f+g}$.
As cohomologies of shifted differentials, they can be assembled to give algebra structures 
on the fibers of $T\M$. From the definitions of $\R_{f+g}$, however, it is not even apparent that
  $\R_{f+g}$ has the same vector space dimension as $\R$.  
  This should be considered a corollary of the theorem \ref{result2} and its proof.
\end{remark}

\section{\v Cech versus Dolbeault}\label{cech}
As mentioned earlier, the Dolbeault resolution 
$(\Gamma(\Theta_X^\bullet\otimes\bar{\Omega}^\bullet_X),\db)$ of
$H^\bullet(\Theta^\bullet_X)$ becomes
a differential graded Lie
algebra with the addition of the Schouten bracket.  
Its cousin $\left(C^\bullet(\Theta^\bullet_X),\dc\right)$, the complex of \v Cech
cochains with values in $\Theta^\bullet_X$, is also a resolution of
$H^\bullet(\Theta^\bullet_X)$ but does not carry a
Lie bracket.  It does, however, have an $L_\infty$ structure,
which, for the purposes of deformation theory, is
just as good.
\subsection{$L_\infty$ algebras}
\subsubsection{Two reminders about vector spaces}

For a graded vector space $V$, we denote by $\sym V$ the graded cofree
cocommutative coalgebra generated by $V$.  That is, $\sym V$ is the
subcoalgebra of the graded tensor coalgebra $\oplus_{k\geq
  0}V^{\otimes k}$, with the standard comultiplication, invarient
under 
the
action of the symmetric group $\Sigma_k$ on $V^{\otimes k}.$

For $i\in \mathbb{Z}$ there is the shift functor $[i]$ acting on a
graded vector space $V=\oplus_k V^k$ defined by $V[i]=\oplus_k V[i]^k$
where $V[i]^k=V^{i+k}$.  The action of the shift functor $[i]$ can be
considered as tensoring with the trivial vector space concentrated in
degree $i.$ 

The exterior and symmetric products are related in this
graded environment by $\sym^i(V[1])\simeq \left(\Lambda^iV
\right)[i].$

\subsubsection{$L_\infty$ algebras}

\begin{definition}An $L_\infty$ algebra consists of a pair $(V,Q)$ where $V$ is
  a graded vector space and $Q$ is a degree one codifferential (square
  zero coderivation) on $\sym (V[1])$.
\end{definition}

Because $\sym (V[1])$ is cofree, any coderivation $Q$ on $\sym (V[1])$
is determined by the pieces (projections followed by restrictions) $Q_i
:\sym^i V[1] \rightarrow V[2-i]$.
The condition that $Q^2=0$ imposes constraints on the linear maps $Q_i.$
The constraints involving the first two components of $Q$ are as
follows: $Q_1:V\rightarrow V[1]$ must be a differential on $V$;
$Q_2:\Lambda^2 V\rightarrow V$ is a skew symmetric bilinear operator
on $V$ and $Q_1$ is a derivation of $Q_2$; and $Q_2$ satisfies the
Jacobi identity up to a homotopy defined by $Q_3:\Lambda^3 V
\rightarrow V[-1]$.  Often, $L_\infty$ algebras are called ``strong
homotopy Lie algebras'' or ``sh Lie algebras.'' The introduction for
physicists in  \cite{sta} is a good introduction for mathematicians,
too.

Any differential graded Lie algebra can be regarded as an $L_\infty$
algebra.  In particular, $\g=\Gamma(\Theta^\bullet_X\otimes
\overline{\Omega}^{\bullet}_{X})$ has an $L_\infty$ structure given by
$Q$ with $Q_1=\db$, $Q_2=[\,,\,],$ and $Q_3=Q_4=\cdots=0.$ However,
even when studying dg Lie algebras, it can be important to view them
as $L_\infty$ algebras when studying maps between them.

\subsubsection{$L_\infty$ morphisms}

\begin{definition}\label{qi} 
An {$L_\infty$ morphism} $(V,Q)\rightarrow (V',Q')$ is a
  differential coalgebra map $\mu:
  (\sym(V[1]),Q) \rightarrow (\sym(V'[1]),Q')$. 
  We say that $\mu$ is a \emph{quasi-isomorphism} if the first
  component of $\mu$, which is a map of complexes
  $\mu_1:(V,Q_1)\rightarrow (V',Q_1')$, induces an isomorphism
  $\mu^*_1:H(V,Q_1)\rightarrow H(V',Q'_1).$ We say that two $L_\infty$
  algebras $(V,Q)$ and $(V',Q')$ are \emph{quasi-isomorphic} if they are
  equivalent under the equivalence relation generated by
  quasi-isomorphisms.
\end{definition}

It is a standard theorem, see for example  \cite{SS, K1}, that if $$\mu:
  (\sym(V[1]),Q) \rightarrow (\sym(V'[1]),Q')$$ is a quasi-isomorphism,
  then there exists a quasi-isomorphism $$\mu':
  (\sym(V'[1]),Q') \rightarrow (\sym(V[1]),Q).$$  
The equivalence relation generated by
  quasi-isomorphisms as stated in the definition of
  \emph{quasi-isomorphic} above is quite natural.

Given an $L_\infty$ algebra $(V,Q)$, one can form the differential
graded Lie algebra, $Coder(\sym(V[1],Q)$, which consists of 
coderivations of $\sym(V[1])$, a differential given by $ad(Q)$, and
the usual Lie
bracket of coderivations.
Then, an $L_\infty$ morphism $(\sym(V[1]),Q) \rightarrow
(\sym(V'[1]),Q')$ is equivalent to a dg Lie morphism
$$Coder(\sym(V[1],Q)\rightarrow Coder(\sym(V'[1],Q').$$  This
description of $L_\infty$ algebras is similar to the ``geometric''
picture promoted in  \cite{K1,Ba,Merk} of $Q$ as an odd, square zero vector
field on a formal pointed graded manifold.

\subsubsection{Deformation theory}\label{deformation_theory}
Given an $L_\infty$ algebra $(V,Q)$, the functor $\mathrm{Def}(V,Q)$
acting from the category of graded local Artin algebras to the
category of sets is defined by
$$\mathrm{Def}(V,Q)(A)=\{\text{coalgebra maps } f:m^*\rightarrow
\sym{V[1]} \text{ with }Q(f(m^*))=0\}/\sim$$
where $m^*$ is the dual of
the maximal ideal of $A$ and $\sim$ is a gauge equivalence, described
geometrically in  \cite{Ba,K1}.  If $(V,Q)$ and
$(V',Q')$ are quasi-isomorphic, then the functors $\mathrm{Def}(V,Q)$
and $\mathrm{Def}(V',Q')$ are canonically equivalent.  In situations
where $\mathrm{Def}(V,Q)$ is representable; that is
$\mathrm{Def}(V,Q)=Hom(\,\cdot\,,\O_\N)$ for some
algebra $\O_\N$, then we call $\N$ the
moduli space.  
The conclusion \cite{SS,K1}  needed for this paper is:
\begin{theorem}
  If $(V,Q)$ and $(V',Q')$ are quasi-isomorphic $L_\infty$ algebras
  and $\mathcal{M}$ is the moduli space for $\mathrm{Def}(V,Q)$ then
  $\mathcal{M}$ is the moduli space for $\mathrm{Def}(V',Q').$
\end{theorem}

We add an opinion here about $Lie$ versus $L_\infty$.  Some
traditional deformation theories are concerned with
$\mathrm{Def}(V,Q)$ in the case that $(V,Q)$ is a differential graded
Lie algebra; that is, $Q_3=Q_4=\cdots=0$.  In this case,
$$\mathrm{Def}(V,Q)(A)=\{\gamma\in (V\otimes m)^1\mid
Q_1(\gamma)+\frac{1}{2}Q_2(\gamma,\gamma)=0\}/\sim.$$
Here,
$\mathrm{Def}(V,Q)$ reduces to the functor described in section \ref{M1}.
But even when $(V,Q)$ is
a dgLie algebra, one cannot effectively study $\mathrm{Def}(V,Q)$
entirely within the more comfortable category of dgLie algebras.  It
is necessary to involve $L_\infty$ structures because $L_\infty$
morphsims provide transformations of (dg Lie) deformation
functors---dg Lie maps are not plentiful enough.  It seems inevitable,
then, that the dg Lie algebra which controls the deformation
theory should be considered an $L_\infty$ object to begin with.

\subsection{$L_\infty$ structure on \v Cech cochains}
In general, the vector space of \v Cech cochains with values in a
sheaf of Lie algebras does not form a Lie algebra, but rather an
$L_\infty$ algebra, canonical up to $L_\infty$ quasi-isomorphism.  The
lemma below explains this phenomenon in some generality.  We will
apply it to the space of \v Cech cochains with values in a holomorphic
vector bundle whose space of sections carries a Lie bracket.

\begin{lemma}\label{export}
  Suppose $(V,Q)$ is an $L_\infty$ algebra and $(V',d)$ is a complex.
  If $(V',d)$ and $(V,Q_1)$ are quasi-isormorphic as complexes, then
  there exists a degree one differential $Q'$ on $\sym V'[1]$ with
  $Q_1'=d$ making $(V,Q)$ and $(V',Q')$ quasi-isormorphic as
  $L_\infty$ algebras.
\end{lemma}

\begin{proof}Here we sketch a simple existence proof.
  If $(V',d)$ and $(V,Q_1)$ are quasi-isomorphic as complexes, then
  the dg Lie algebras $$(Coder(\sym V'[1]),\mathrm{ad}(d))\text{ and
    }(Coder(\sym V[1]),\mathrm{ad}(Q_1))$$
  will be quasi-isormorphic
  as dgLie algebras (one must use the cofreeness of the symmetric product).  
    Since $\mathrm{ad}(Q_1)$ extends to an
  integrable differential $\mathrm{ad}(Q)$ on $(Coder(\sym V[1])$, so
  must $\mathrm{ad}(d)$ extend to a differential
  $D'_{Q'}=\mathrm{ad}(Q')$ making the dg Lie algebras
  $$(Coder(\sym V'[1]),\mathrm{ad}(Q'))\text{ and }(Coder(\sym
  V[1]),\mathrm{ad}(Q))$$
  quasi-isomorphic.
\end{proof}

Suppose that $A$ is a holomorphic vector bundle with
$(\Gamma(A),[\,,\,])$ a Lie algebra.  The bracket on $\Gamma(A)$ can
be extended using the associative $\wedge$ product on
$\Gamma(\bar{\Omega}^\bullet)$ to a Lie bracket on $\Gamma(A\otimes
\bar{\Omega}^\bullet)$.  Since $(A\otimes \bar{\Omega}^\bullet,\db)$
is a fine resolution of $A$, the \v Cech complex
$\left(C^\bullet(A),\dc\right)$ and the Dolbeault complex
$(\Gamma(A\otimes \bar{\Omega}^\bullet),\db)$ are quasi-isomorphic.
Then, Lemma \ref{export} can be used to export the Lie algebra
structure from the Dolbeault complex to an $L_\infty$ structure (which
will not in general be Lie), to the \v Cech complex.

Applying this when $A=\Theta_X^\bullet$, we obtain the following
theorem.  Recall that $\g$ denotes $\Gamma\left(\Theta_X^\bullet
  \otimes \bar{\Omega}_X^\bullet\right)$.

\begin{theorem}\label{thm6}
  There exists an $L_\infty$ structure $\check{Q}$ on
  $C^\bullet(\Theta^\bullet_X)$, with $\check{Q}_1=\dc$, so that
  $$(C^\bullet(\Theta^\bullet_X),\check{Q}) \text{ and
    }(\g,\db,[\,,\,])$$
  are quasi-isomorphic.
\end{theorem}

\subsubsection{The $\check{Q}_2$ term}
It is possible to give a constructive proof of Lemma \ref{export}.
That is, $d$ can be extended to $Q'=d+Q'_2+Q'_3+\cdots$
perturbatively.  See, for example, the articles \cite{BFLS, Merk} for
other explicit, perturbative constructions of $L_\infty$ structures.
Here is an outline for producing $Q'_2$.

For simplicity, and because this is the case of interest, let
$(V,Q=Q_1+Q_2)$ be a differential graded Lie algebra and suppose that
$\phi:(V',d)\rightarrow (V,Q_1)$ is a quasi-isormorphism of complexes.
Then, there exists a quasi-isomorphism $\psi:(V,Q_1)\rightarrow
(V',d)$ which is a homotopy inverse of $\phi$.  In particular, there
exists an $s:V \rightarrow V[-1]$ with
$$\phi\psi -1=sQ_1+Q_1s.$$
Now, one can define a skew-symmetric
bilinear map $Q'_2:V'\otimes V'\rightarrow V'$ by
$$Q'_2(v_1,v_2)=\psi (Q_2(\phi(v_1),\phi(v_2))).$$
Because $\psi$ and
$\phi$ are inverses only up to homotopy, $Q'_2$ will not satisfy the
Jacobi identity.  However $Q'_2$ will satisfy Jacobi up to homotopy,
which is all that is required for the second term of and $L_\infty$
structure on $V'$.  The higher homotopies are constructed inductively.
Also, note that $\phi$ is not a homomorphism in the sense that
$\phi(Q'_2(v_1,v_2))\neq Q_2(\phi(v_1),\phi(v_2))$; but $\phi$ does
define an $L_\infty$ morphism, which can be thought of as a
homomorphism up to homotopy.  In fact, the collection of maps
$F_k:\Lambda^k V' \rightarrow V$ given by
\begin{align*}
  F_1&=\phi \\
  F_2&= s \circ Q_2 \circ (\phi \otimes \phi)\\
  \vdots &
\end{align*}
assemble to define an $L_\infty$ morphism from $(V',d+Q'_2+\cdots)$ to
$(V,Q_1+Q_2).$

This outline describes $Q'_2$ modulo the precise information about the
quasi-isomorphisms between $(V,Q_1)$ and $(V',d)$.  For the example of
\v Cech cochains, $\check{Q}_2$ can be described explicitly once the
quasi-isomorphisms between $\left(C^\bullet(A),\dc\right)$ and
$(\Gamma(A\otimes \bar{\Omega}^\bullet),\db)$ are analyzed.

\begin{remark*}
  The construction of $Q'_2$ above depends on the choice of $\psi$ and
  $s$.  However, a change of choice induces a quasi-isomorphism
  between the respective $L_\infty$ structures.
\end{remark*}

 \subsubsection{The quasi-isomorphisms between \v Cech and Dolbeault}
 \label{quisms}
 
 As before, assume that $A$ is a holomorphic vector bundle and that
 $(\Gamma(A),[\,,\,])$ is a Lie algebra.  Consider the resolution of
 $A$ by the complex of fine sheaves $(A^\bullet,\db)$ where
 $A^p=A\otimes \overline{\Omega}^{p}_X$.  The inclusions
 $$
 \xymatrix{ & \ar@{<-}[dl] (C^\bullet(A^\bullet),D=\dc+\db)
   \ar@{<-}[dr] \\
   (C^\bullet(A),\dc) & & (\Gamma(A^\bullet),\db) }
 $$
 are quasi-isormophisms \cite{GH}.  Each map is invertible, up to
 homotopy, so we have maps
\begin{gather*}
  \phi:(C^\bullet(A),\dc) \rightarrow (\Gamma(A^\bullet),\db)\\
  \psi:(H^0(A^\bullet),\db) \rightarrow (C^\bullet(A),\dc)
\end{gather*}
inducing (inverse) isomorphisms in cohomology.

The map $\phi$ is straightforward to describe on the subspace
$Z^\bullet(A).$ Let $a \in C^p(A)$ with $\dc(a)=0$.  One can find
(using a partition of unity) an $\alpha_{p-1}\in C^{p-1}(A^0)$ with
$\dc(\alpha_{p-1})=-a.$ Then,
$$a+D(\alpha_{p-1})=\db(\alpha_{p-1})\in C^{p-1}(A^1).$$
Now,
$\db(\alpha_{p-1})$ is $\dc$ closed since
$$\dc(\db(\alpha_{p-1}))=\db\dc(\alpha_{p-1})=-\db(a)=0.$$

Again, (using a partition of unity), there exists an $\alpha_{p-2}\in
C^{p-2}(A^1)$ with $\dc(\alpha_{p-2})=-\db(\alpha_{p-1}).$ One gets
$$a+D(\alpha_{p-1}+\alpha_{p-2})=\db(\alpha_{p-2})\in
C^{p-2}(A^2)\text{ with }\dc\db(\alpha_{p-2})=0.$$

Continuing, one arrives at $\alpha_0 \in C^0(A^{p-1})$, with
$$a+D(\alpha_{p-1}+\alpha_{p-2}+\cdots +\alpha_0)=\db(\alpha_0) \text{
  and }\dc\db\alpha_0=0.$$
Since $\dc(\db(\alpha_0))=0$,
$\db(\alpha_0)\in C^0(A^p)$ gives a global section $\alpha\in
\Gamma(A^p).$ This describes the map $\phi$:
$$\phi(a)=\alpha.$$

The map $\psi$ is defined on the cycles in $\Gamma(A^\bullet)$ in a
similar fashion.  If $\beta\in \Gamma(A^p)$ with $\db(\beta)=0$, then
one can find (using the Poincar\'e lemma) elements
$$b_{p-k}\in C^{k-1}(A^{p-k})\text{ with
  }\db(b_{p-k})=-\dc(b_{p-k+1}).$$
Then
$$\beta+D(b_{p-1}+b_{p-2}+\cdots + b_0)=\dc(b_0)\in C^p(A^0).$$
Since
$\db\dc b_0=0$, $\dc b_0$ can be identified with an element $b\in
C^p(A).$ We have $$\psi(\beta)=b.$$

\subsubsection{Computing $\check{Q}_2$}

Now, $\check{Q}_2:Z^\bullet(A)\otimes
Z^\bullet(A^\bullet)\rightarrow Z^\bullet(A\bullet)$ is defined by
$$\check{Q}_2(a,b)=\psi[\phi(a),\phi(b)]$$ where the bracket on the
right is the
Schouten bracket in $\Gamma(A^\bullet).$  

\begin{lemma}\label{eq2} 
For $a=\{a^i\}\in Z^0(A^\bullet)$ and $b=\{b^{i_0\cdots i_p}\}\in Z^p(A^\bullet)$
\begin{equation*}
\check{Q}_2(a,b)=c=\{c^{i_0\cdots i_p}\}\in Z^p(A)\text{ where }c^{i_0\cdots
  i_p}=[a^{i_0},b^{i_0\cdots i_p}].
\end{equation*}
\end{lemma}
\begin{proof}This is a computation.
We introduce as a notational aid, 
the bilinear function
$[|\,,\,|]:C^p(A^\bullet)\times C^q(A^\bullet)\rightarrow C^{p+q}(A^\bullet),$ defined for $\alpha\in C^p(A^\bullet)$ and $\beta\in C^q(A^\bullet)$ by 
$$[|\alpha,\beta|]=\gamma=\{\gamma^{i_0\cdots i_p \cdots i_{p+q}}\} \text{ where }
\gamma^{i_o\cdots i_p \cdots i_{p+q}}=[\alpha^{i_0\cdots i_p},\beta^{i_p\cdots i_{p+q}}].$$  
While $[|\,,\,|]$ does not satisfy the
Jacobi identity, nor the Jacobi identity up to homotopy (it is not even skew-symmetric for arbitrary elements of
$C^p(A^\bullet)\otimes C^p(A^\bullet)$.), it does behave well
with respect to the differentials---both differentials are derivations
of $[|\,,\,|]$:
\begin{gather*}
\dc[|a,b|]=[|\dc(a),b|]+(-1)^{|b|}[a,\dc(b)|]\\
\db[|a,b|]=[|\db(a),b|]+(-1)^{|b|}[a,\db(b)|].
\end{gather*}
In this notation, the lemma asserts that for the special case of $a\in
Z^0(A^\bullet)$ and $b\in Z^p(A^\bullet)$, we
have
$\check{Q}_2(a,b)=[|a,b|].$

Let $a=\{a^i\}\in Z^0(A)$ and $b=\{b^{i_0\cdots i_p}\}\in
Z^p(A)$.  Select elements $\beta_0,\beta_1,\ldots,\beta_{p-1}$ with
$$\beta_k \in C^k(A^{p-k-1}),\, \dc \beta_{k-1}=\db \beta_k,\,
\dc\beta_{p-1}=b,\beta=\db \beta_0.$$
We have $\phi(b)=\beta \in \Gamma(A^p)$
and $\phi(a)=\alpha\in \Gamma(A^0)$
identified with $a\in C^0(A)$.  Then $[\alpha,\db\beta_0]\in
\Gamma(A^p)$ can be considered as the cocylce $[|\alpha,\db\beta_0|]\in
C^0(A^p).$ 
Now, we compute $c=Q'_2(a,b)=\psi([\alpha,\beta_0])$.  We set $\gamma
=[\alpha,\db\beta_0]\in \Gamma(A^p)$ and find elements
$\gamma_{0},\gamma_1,\ldots \gamma_{p-1}$
satisfying $$-\dc \gamma_k =\db \gamma_{k+1}\text{ and } \db \gamma_0=\gamma.$$
Then $-\dc \gamma_{p-1}=\psi(\gamma)=c\in C^p(A).$  One finds that 
\begin{align*}
\gamma&=[\alpha,\db\beta_0]\\
&=\db[|a,\beta_0|]\\
&=\db \gamma_0 \text{ for }\quad \gamma_0=[|a,\beta_0|].
\intertext{then}
\dc \gamma_0&=\dc [|a,\beta_0|]\\
&=[|a,\dc\beta_0|]\\
&=-[|a,\db\beta_1|]\\
&=-\db[|a,\beta_1|]\\
&=\db \gamma_1,\text{ for } \gamma_1=[|a,\beta_1|].
\intertext{ and so on.}
\end{align*}
One finds that $\gamma_k=[|\alpha,\beta_k|]$ and at the end 
$$-\dc \gamma_{p-1}=\dc [|a,\beta_{p-1}|]=[|a,\dc
\beta_{p-1}|]=[|a,b|]$$ as claimed.
\end{proof}

It is also straightforward, but not necessary for our purposes, to calculate $\check{Q}_2(a,b)$ for $a\in Z^p(A)$,
$p>0$.  One finds that $\check{Q}_2$ can be expressed in terms of
$[|\,,\,|]$, but not with just \emph{one} such term:  
For example, if $\alpha,\beta \in Z^1(A)$,
\begin{gather}\check{Q}_2(\alpha,\beta)=\frac{1}{2}\left([|\alpha,\beta|]+[|\beta,\alpha|]\right)\label{eq3}\\
\intertext{or, equivalently,}
\check{Q}_2(\alpha,\beta)=\gamma \in Z^2(A) \quad \text{where }\gamma^{ijk}=\frac{1}{2}\left([\alpha^{ij},\beta^{jk}]+[\beta^{ij},\alpha^{jk}]\right).
\end{gather}
The reader may recognize this as classical bracket in the Kodaira-Spencer
theory of deformations of complex structure \cite{Kod}.

\subsection{Associative algebra structures}
If the space of sections of $A$ carries an associative product, then
both the Dolbeaut and \v Cech complexes become associative algebras
and the maps $\phi$ and $\psi$ are algebra maps.

Since $\Gamma\left(A\otimes \bar{\Omega}^\bullet\right)$ is simply the
tensor product of two associative algebras ($\bar{\Omega}^\bullet$ has
the wedge product), it is an associative algebra.  The standard way to
make the \v Cech cochains $C^\bullet(A)$ into an associative algebra
is as follows:
\begin{gather*}
  \text{For }a=\{a^{i_0 \cdots i_p}\}\in C^p(A)\text{ and }b=\{b^{i_0
    \cdots i_p}\}\in C^q(A) \intertext{define} a\cdot b=c\in
  C^{p+q}(A), \text{ where }c^{i_0 \cdots i_{p+q-1}}=a^{i_0 \cdots
    i_p}\cdot b^{i_p \cdots i_{p+q}}.
\end{gather*}
Both differentials are derivations of these algebra structures; so
$\Gamma(A\otimes \bar{\Omega}^\bullet, \db, \wedge)$ and
$(C^\bullet(A), \dc, \cdot)$ are quasi-isomorphic as differential
graded associative algebras.

In the case that the sections of $A$ have both an associative product
and a Lie bracket satisfying the Poisson identity equation
\ref{gerst}, then the Dolbeaut resolution $A\otimes
\bar{\Omega}^\bullet$ will also be a Gerstenhaber algebra.  The
$\check{Q}_2$ term of the $L_\infty$ structure on $C^\bullet (A)$,
being a sum of terms of the form $|[\,,\,|]$, also satisfies a Poisson
identity.

One way that $C^\bullet (A)$ becomes an algebra is when
$A=\Lambda^\bullet B$ for some bundle $B$.  Then the wedge product on
sections of $A$ induces a wedge product on $C^\bullet (A)$, which is
usually also denoted by $\wedge$.  This is the case for
$A=\Theta^\bullet_X$.  A common example occurs for $A=\Lambda^\bullet
\Omega_X$ and the maps $\phi$ and $\psi$ induce \emph{multiplicative}
isomorphisms between $H^\bullet(\Omega_X^\bullet)$ and
$H^{\bullet,\bullet}(X).$

\section{A model for $\Theta_X$}\label{L}

In this section, we construct a sheaf of differential graded Lie
algebras $\L$ on $\P^n$ that is quasi-isomorphic to $\Theta_X$ (with
the zero differential and the ordinary bracket of vector fields).

Let $T=S[y]=S[x^0,\ldots, x^n,y]$ and with weights and degrees as follows: 
\begin{align*}
\w(x)&=1 &  \w(y)&=\nu \\
\deg(x)&=0 & \deg(y)&=-1.
\end{align*}
As before, we assume that $f\in S$ has been chosen, that
$\w(f)=\nu=n+1$, and that $f$ is nonsingular.
Consider $\der(T)$.  We use the notation $\partial_i$
for $\frac{\partial}{\partial x^i}$, $\partial$ for
$\frac{\partial}{\partial y}$, and $g_i$ for $\partial_i g$.
We have weights and degrees:
\begin{align*}
\w(\partial_i)&=-1 &  \w(\partial)&=-\nu \\
\deg(\partial_x)&=0 & \deg(\partial)&=1.
\end{align*}
We define\footnote{One way to view $d_f$ is to define a differential in $T$ by $d(y)=f$
and $d(x^i)=0$.  Then $d_f=\mathrm{ad}(d)$ in $\der(T).$} an $S$ linear, weight zero, degree one, differential
  $d_f:\der(T)\rightarrow \der(T)$
  by
\begin{align*}
d_f(y\partial_i)&=f\partial_i-f_i y \partial \\
d_f(\partial_i) &= f_i \partial \\
d_f(y \partial) &= f\partial.
\end{align*}
The triple $(\der(T),d_f,[\,,\,])$ is a differential graded Lie
algebra, where the grading is by \emph{degree} and $[\,,\,]$ is the ordinary
bracket of derivations.  That $\der(T)$ is weighted is a feature
that will be useful for sheaf computations later.

Now, adjoin a weight zero, degree zero element $e$ with
\begin{gather*}
d_f(e)=\sum_{k=0}^n x^k \partial_k -\nu y\partial \\
[e,a]=[a,e]=\w(a) \text{ for any }a\in \der(T)\oplus Se.
\end{gather*}
The result $(L,d_f)$ where $L=\der(T)\oplus Se$ is a graded Lie algebra and a
complex (and the
weight zero subcomplex $\{a \in \der(T)\oplus
Se \mid \w(a)=0\}$ is a differential graded Lie algebra when
considered with the bracket).

\subsection{Cohomologies of $L$}

Here is a chart to help keep track of degrees and weights:

\medskip
\begin{center}
\begin{tabular}{|c|c|c|c|c|c|}
\hline 
& $e$ & $y\partial_i$ & $\partial_i$ & $y\partial$ & $\partial$ \\
\hline 
degree & $-1$ & $-1$ & $0$ & $0$ & $1$   \\ 
\hline 
weight & $0$ & $\nu -1$ & $-1$ & $0$ & $-\nu$ \\
\hline 
$d_f$ & $\ds \sum_{k=0}^n x^k \partial_k -\nu y \partial$ &
$f\partial_i-f_iy\partial$ & $f_i\partial$ & $f\partial$ & $0$ \\
\hline
\end{tabular}
\end{center}
\medskip
One reads $H^{-1}(L)=0$ and $H^1(L)=S\partial/J_f\partial=S/J_f$ and
$$H^{0}(L)=\left\{\sum_{i=0}^n g^i \partial_i \mid \sum_{i=0}^n g^i
  f_i =hf \right\}\biggr{/}\left\{ p \left(\sum_{i=0}^n
    x^i\partial_i\right)\right\}$$
where $g^i,h,p\in S$.

Let $\L$ denote the sheaf of differential graded Lie algebras on
$\P^n$ associated to $L$:
$$\Gamma(\L,U_i)=\left\{\frac{\alpha}{(x^i)^r}: \alpha\in L,
  \w(\alpha)=r\Leftrightarrow \w\left(\dfrac{\alpha}{(x^i)^r}\right)=0
\right\}.$$  Note that $\H^0(\L)\simeq \Theta_X$.  Here caligraphic
$\H$ 
stands for the kernel of $d_f$ modulo the image of $d_f$ 
(not the hypercohomology of $\L$).  Since $H^1(L)$ is finite dimensional over $\mathbb{C}$, for
the complex of sheaves $\L$, we have $\H^1(\L)=0$.  So
$\H^{-1}(\L)=\H^1(\L)=0$ giving $\H(\L)=\H^0(\L).$ That is,
$$\H(\L)\simeq \Theta_X.$$

\subsection{A simple formality for $\L$}
For any three term sheaf of dgLa's 
$$0\rightarrow \A^{-1}\overset{d^0}{\rightarrow} \A^0
\overset{d^1}{\rightarrow} \A^1\rightarrow 0$$  with $\H^{-1}(\A)=\H^1(\A)=0$, 
the vertical maps $\eta$ and $\zeta$
$$
\xymatrix{ \A^{-1} \ar@{->}[r]^{d^0} & \A^0
 \ar@{->}[r]^{d^1}&\A^1\\
\ar@{->}[u]_\zeta\ar@{->}[d]_\eta \A^{-1} \ar@{->}[r]^{d^0} & \ar@{->}[u]_\zeta\ar@{->}[d]_\eta  \ker(d^1)
 \ar@{->}[r]^{d^1}&\ar@{->}[u]_\zeta \ar@{->}[d]_\eta 0\\
0\ar@{->}[r]^{d^0} & \H(\A)
\ar@{->}[r]^{d^1}&0}
$$
are
quasi-isomorphisms of sheaves of dgLa's (meaning the induced map
on $\H$ is an isomorphism). This is the situation for $\L$.

\begin{theorem}\label{thetathm}
 $(\L,d_f,[\,,\,])$ is quasi-isomorphic as a sheaf of dgLa's to $(\Theta_X,0,[\,,\,])$.
\end{theorem}

\section{The complex $(\C,D)$ and the $L_\infty$ algebra $(\C,Q^{\C})$}\label{F}\label{theta}
Let $$F=\bigoplus_{k=0}^{n-1} F^k,\quad F^k=\sym^k(L[1])\simeq
\left(\Lambda^k L\right)[k].$$
The differential on $L$ extends (by the
Leibnitz rule) to $F$ and the bracket on $L$ extends to a
Schouten--type bracket on $F$.  They are determined by
\begin{itemize}
\item for $\alpha\in F^i$, $\beta\in F^j$, $d_f(\alpha \wedge
  \beta)=d_f(\alpha)\wedge \beta + (-1)^i \alpha \wedge d_f(\beta)$,
\item for $g\in F^0=T=S[x^0,\ldots, x^n,y]$ and $\alpha \in F^1=L[1],$
  $[\alpha,g]=\alpha(g)$,
\item for $\alpha\in F^i$, $\beta\in F^j$, $\gamma\in F$, $[\alpha
  \wedge \beta,\gamma]=\alpha \wedge [\beta,\gamma]+(-1)^{(i+1)j}
  [\alpha,\gamma]\wedge \beta.$
\end{itemize}

Again, here is a chart of degrees and weights.  Note that $x^i,y \in
F^0\subset F$ and remember that $L$ was shifted: \medskip
\begin{center}
\begin{tabular}{|l|c|c|c|c|c|c|c|c|c|}
\hline 
 & $y$ & $x^i$ & $e$ & $y\partial_i$ & $\partial_i$ & $y\partial$ &
$\partial$ &  $\partial^l $ \\
\hline 
deg & $-1$ & $0$ & $0$ & $0$ & $1$ & $1$ & $2$ & $2l$ \\
\hline 
w & $\nu$ & $1$ & $0$ & $\nu-1$ & $-1$ & $0$ & $-\nu$ &
 $-l\nu$ \\
\hline 
\end{tabular}
\end{center}
\medskip

Let $\F$ denote the sheaf associated to $F$ and denote by $C^p(\F)$,
the \v Cech $p$-cochains with values in $\F.$ Let
$\dc:C^p(\F)\rightarrow C^{p+1}(\F)$ denote the \v Cech differential.
We have a double complex whose $(p,q)$ term is $C^p(\F^q)$ and a
single complex $(\C,D)$ where $\C^k=\oplus_{p+q=k}C^p(\F^q)$ and
$D=\dc+d_f.$

\begin{theorem}\label{thm7}
  The complexes $(\C,D)$ and $\left(\g,\db\right)$ are
  quasi-isomorphic.
\end{theorem}

\begin{proof}
  From theorem \ref{thetathm}, the complexes of sheaves $(\L,d_f)$ and
  $(\Theta_X,0)$ are quasi-ismorphic.  Therefore, by taking exterior
  powers, $(\F,d_f)$ and $(\Theta_X^\bullet,0)$ are quasi-ismorphic
  (as complexes of sheaves).  Since quasi-isomorphic complexes induces
  isomorphisms in hypercohomology \cite{GH}, $(C^\bullet(
  \F),D=\dc+d_f)$ and $(C^\bullet(\Theta_X^\bullet),\dc+0)$ are
  quasi-isormorphic as complexes.
  
  By theorem \ref{thm6}, we have a quasi-isomorphism between
  $(C^\bullet(\Theta_X^\bullet),\dc)$ and \linebreak
  $\left(\Gamma\left(\Theta^p_X\otimes
      \overline{\Omega}^q_X\right),\db\right)$, establishing the
  theorem.
\end{proof}

There are advantages to working with $(\C,D)$.  For one, we are able 
to compute the cohomology of $(\C,D)$.  
Another is that $\C$ is defined independently of
$X$---the polynomial $f$, and hence all of the geometry of $X$, is
carried in the differential $D$.  When the differential $D$ is
deformed, it is not difficult to interpret the shifted cohomology
rings.

\begin{theorem}\label{thm8}
  $H(\C,D)\simeq \R$.
\end{theorem}

\begin{proof}
  Consider the double complex $(C^\bullet(\F^\bullet),\dc,d_f)$ and
  the filtration on the single complex $(\C^\bullet,D):$
  $$F^p\C^\bullet= C^0(\F^{\bullet})\oplus C^1(\F^{\bullet-1})\oplus
  \cdots \oplus C^{\bullet-p}(\F^{p}).$$
  Let $\{E_r,\delta_r\}$ denote
  the spectral sequence associated to this filtered complex.  We will
  see that spectral sequence degenerates at the $E_2$ term:
  $E_\infty^{p,q}=E_2^{p,q}\simeq H^q_{d_f}(H^p_{\dc} (\F)).$
\begin{figure}[h]
\begin{center}
  \setlength{\unitlength}{0.0017cm}
\begin{picture}(5500,2500)(-100,0)

\put(500,1800){\blacken\ellipse{100}{100}}
\put(1000,200){\blacken\ellipse{100}{100}}
\put(2000,200){\blacken\ellipse{100}{100}}
\put(3000,200){\blacken\ellipse{100}{100}}
\put(4000,200){\blacken\ellipse{100}{100}}

\path(500,0)(500,800)\dottedline{75}(500,800)(500,1200)\path(500,1200)(500,2100)
\path(1000,0)(1000,800)\dottedline{75}(1000,800)(1000,1200)\path(1000,1200)(1000,2100)
\path(2000,0)(2000,800)\dottedline{75}(2000,800)(2000,1200)\path(2000,1200)(2000,2100)
\path(3000,0)(3000,800)\dottedline{75}(3000,800)(3000,1200)\path(3000,1200)(3000,2100)
\path(4000,0)(4000,800)\dottedline{75}(4000,800)(4000,1200)\path(4000,1200)(4000,2100)

\path(200,200)(4300,200)\dottedline{75}(4300,200)(4600,200)
\path(200,500)(4300,500)\dottedline{75}(4300,500)(4600,500)
\path(200,1500)(4300,1500)\dottedline{75}(4300,1500)(4600,1500)
\path(200,1800)(4300,1800)\dottedline{75}(4300,1800)(4600,1800)

{\footnotesize
\put(4900,-50){\makebox(0,0)[lm]{$p$ ($=$ deg in $\F$)}}
\put(100,200){\makebox(0,0)[rm]{$0$}}
\put(100,500){\makebox(0,0)[rm]{$1$}}
\put(100,1500){\makebox(0,0)[rm]{$n-1$}}
\put(100,1800){\makebox(0,0)[rm]{$n$}}

\put(100,2300){\makebox(0,0)[cb]{\hbox{$q$ ($=$ \v Cech deg)}}}
\put(500,-50){\makebox(0,0)[ct]{$-1$}}
\put(1000,-50){\makebox(0,0)[ct]{$0$}}
\put(2000,-50){\makebox(0,0)[ct]{$2$}}
\put(3000,-50){\makebox(0,0)[ct]{$4$}}
\put(4000,-50){\makebox(0,0)[ct]{$6$}}
\put(2100,300){\makebox(0,0)[lb]{$a$}}
\put(600,1750){\makebox(0,0)[lt]{$e_k$}}
}
\end{picture}
\end{center}
\caption{$E_1^{p,q}$}\label{fig2}
\end{figure}

We compute $E_1^{p,q} \simeq H^q(\F^p,\dc)$.  To describe the classes
in $E_1=\oplus_q H^q_{\dc}(\F)$ that survive in $E_2$ (see figure 1),
first note that since the $\F^p$ are locally free sheaves on $\P^n$,
$$H_{\dc}^q(\F)=0, \text{ unless $q=0$ or $n$.}$$

\paragraph{Case $q=0$.}    The (homogeneous) elements generating $H^0(\F)$ that
are $d_f$ closed and not necessarily $d_f$ exact are
given by $a=\{a^i\}\in C^0(\mathcal{U},\F)$ where
\begin{equation*}
a^i=\dfrac{x^i g(x)\partial^p}{x^i} \in \Gamma(U_i,\F), \quad
g(x)\in S^{p\nu}.
\end{equation*} 
An $a$ as above
is a $d_f$ boundary if and only if $g(x)=\sum_{k=0}^n g^i(x)f_i(x)$
for some $g^i\in S^{(n-1)\nu+1}.$  So,
$$H^{2p}_{d_f}(H^0_{\dc}(\F))\simeq R^{p}=\im\left
  (S^{[p\nu]}\rightarrow S/<f_1,\ldots,f_n>\right)$$
and $$\oplus_{p} H^{2p}_{d_f}(H^0_{\dc}(\F))\simeq \oplus_p R^p.$$

\paragraph{Case $q=n$.} Each class in $H_{\dc}^{n}(\F)$ is of the form $$e_k\overset{\mathrm{def}}{=}\frac{ye^k}{x^0 \cdots
  x^n},\quad k=0,1,\ldots , n-1.$$
Now,
$$d_f\left(\frac{ye^k}{x^0\cdots x^n}\right)=\frac{fe^k - y
ke^{k-1}\left(\sum_{i=0}^n x^i \partial_i\right)}{x^0 \cdots x^n}.$$
Because $d_f(e_k)$ can be expressed as a sum of terms, each of which has a
positive power of an $x^i$ in the numerator,
$$d_f(e_k)\in \dc \left(
  C^{n-1}(\F)\right).$$ 
Therefore, each $e_k$ is a class in $E_2$ of total degree $n-1$ (\v Cech degree $n$ and degree $-1$ in
$\F$)
and
$$\oplus H^{-1}_{d_f}(H^n_{\dc}(\F))\simeq \mathbb{C} e_k.$$ 
Note, too, that it is impossible for $D(e_k)$ to kill any of the \v
Cech degree zero classes, hence $E_2=E_3=E_4=\cdots.$

Together, we have
$$H(\C,D)\simeq \oplus_{p} H^{2p}_{d_f}(H^0_{\dc}(\F))\oplus
H^{-1}_{d_f}(H^n_{\dc}(\F))\simeq \oplus_p R^p\oplus_{k=0}^{n-1} \mathbb{C} e_k=\R.$$ 
\end{proof}

\subsubsection{Moduli}

Since, $(\C,D)$ is quasi-isomorphic to $(C(\Theta_X^\bullet,\dc))$,
there exists an $L_\infty$ structure, call it $Q^{\C}:\sym \C[1]\rightarrow \sym
\C[1]$, with $Q^\C_1=D$ such that 
$(\C,Q^\C)$ is quasi-ismorphic to $(C(\Theta_X^\bullet,\check{Q})$. 
Since $(C(\Theta_X^\bullet,\check{Q})$ is quasi-ismorphic to
$(\g,\db,[\;,\;])$, we have

\begin{theorem}
As $L_\infty$ algebras, 
$\ds( \C,Q^{\C})$ and $(\g,\db,[\,,\,])$ are quasi-isomorphic.
\end{theorem}
\begin{corollary}
The moduli space for $\mathrm{Def}(
\g,\db,[\,,\,])$ equals the moduli space for
$\mathrm{Def}(\C,Q^{C})$.
\end{corollary}

This completes the proof of theorem $\ref{result1}$.

In order to complete the proof of theorem \ref{result2}, 
we analyze
the cohomology of $H(\C,D_b)$ for the shifted differentials $D_b$.

\subsubsection{Shifting $D$}
We
take $[a] \in H(C,D)\simeq \R$ and obtain
a shifted $L_\infty$ algebra $(\C,Q^{\C,a})$.  The linear part of
$Q^{\C,a}$ is a shifted differential \cite{Merk} $D_a:=Q_1^{\C,a}$
defined by
$$Q_1^{\C,a}=D+Q_2^\C(a,\cdot)+\frac{1}{2!}Q_3^{\C}(a,a,\cdot)+\frac{1}{3!}Q_4^{\C}(a,a,a,\cdot)+\cdots$$
The classes in $H(\C,D)\simeq \R$ come in two types: 
\begin{enumerate}
\item the primite
elements $[a]$,
\begin{equation}\label{h0part}
a^i=\dfrac{x^i g(x)\partial^p}{x^i} \in \Gamma(U_i,\F), \quad
g(x)\in S^{p\nu}
\end{equation} 
have \v Cech degree zero and even degree $2p$ in $\F$
\item
the nonprimitive elements $[e_k]$ have \v Cech degree $n$ and degree
$-1$ in $\F$.  
\end{enumerate}
For a primitive $a$, a simple calculation shows that
$$Q_2^\C(Q_2^\C(a,a),b)+Q_2^\C(Q_2^\C(b,a),a)+Q_2^\C(Q_2^\C(a,b),a)=0.$$
That is, $Q_2^\C$ satisfies the Jacobi identity (exactly, not just up
to homotopy) when two of the three elements are \v Cech zero cocycles
in $\C$.  This implies that $Q$ may be taken with
$$Q_3^\C(a,a,b)=Q_4^\C(a,a,a,b)=\cdots =0.$$
So, $D_a$ reduces to the familiar
$$D_a=D+Q^{\C}_2(a,\cdot).$$
For the nonprimitive elements $e_k$,
degree restrictions imply that
$$Q_3^\C(e_k,e_k,b)=Q_4^\C(e_k,e_k,e_k,b)=\cdots =0,$$
and again we
have
$$D_{e_k}=D+Q^{\C}_2(e_k,\cdot).$$
Now, we recompute the cohomologies
$H(\C,D_a)$ and $H(\C,D_{e_k})$:

Consider a $(0,2p)$ cocycle $a\in \C$ and let $g\in S^{[p\nu]}$ as in
Equation \ref{h0part}.  For $b\in \C$, we can determine $Q_2(a,b)$ by
$[a^{i_0},b^{i_0\cdots i_q}]$ where this bracket is just the bracket
in $\F$.  Since $\Gamma(\mathcal{U},\F)$ is generated by symmetric
products of elements of $\Gamma(\mathcal{U},\F^i)$, $i=0,1$, (and
$[|\,,\,|]$ satisfies the Liebnitz property with respect to these
symmetric products) we need only compute $[a,b]$ for generators $b$.
The results
\begin{align*}
  [g \partial^p, y \partial_i] &= \left(g\partial_i - \frac{\partial
      g}{\partial x^i}y\partial\right) \partial^{p-1}\\
  [g \partial^p, \partial_i] &= \left(\frac{\partial g}{\partial
      x^i}\partial\right)\partial^{p-1}\\
  [g \partial^p, y \partial] &= \left(g \partial\right)\partial^{p-1}\\
  [g \partial^p, \partial] &= 0\\
  [g \partial^p, e] &= 0.
\end{align*}
yield
\begin{lemma}\label{shifta}For $a$ as in Equation \ref{h0part}
  $$D_a:=\dc+d_f+Q^{\C}_2(g,\cdot)=\dc+d_{f+g}$$
  where
  $d_{f+g}(b)=d_f(b)+d_g(b)\partial^{p-1}.$
\end{lemma}

Lemma \ref{shifta} above and the same specral sequence computation
(cf., the proof of theorem \ref{thm8}) that gives $H(\C,D)\simeq \R$,
except with $f$ replaced by $f+g\partial^p$, proves that
$H(\C,D_a)\simeq \R_{f+g}$.

To determine $Q_2(e_k,b)$, we need only look at $Q_2(e_k,b)$ for
$b\in H^0_{\dc}
(\F)$, since $Q_2(e_k,b)\in \im(D) $ if the \v Cech degree of $b>0$:
\begin{align*}
  \left[\frac{ye^k}{x^0 \cdots x^n}, y \partial_i\right] &= 0\\
  \left[\frac{ye^k}{x^0 \cdots x^n}, x^j\partial_i\right] &=
  \frac{x^jye^k}{x^0\cdots x^n\cdot x^i}\in \dc \left(C^{n-1}(\F)\right)\\
  \left[\frac{ye^k}{x^0 \cdots x^n}, y \partial\right] &=
  -\frac{ye^k}{x^0 \cdots x^n}\\
  \left[\frac{ye^k}{x^0 \cdots x^n}, g \partial\right] &= -\frac{ge^k}{x^0 \cdots x^n}\in \dc \left(C^{n-1}(\F)\right)\\
  \left[\frac{ye^k}{x^0 \cdots x^n}, e\right] &= 0.
\end{align*}
Therefore, $D$ closed elements are $D_{e^k}$ closed and $D$ exact 
elements are $D_{e^k}$ exact elements:
$$H(\C,D_{e_k})=H_{d_f+\mathrm{ad}(e_k)}H^\bullet_{\dc}
(\F)=H_{d_f}H^\bullet_{\dc} (\F)=H(\C,D)=\R.$$
Together,
\begin{theorem}
  $H(\C,D_{e_k})\simeq \R$ and $H(\C,D_{a})\simeq \R_{f+g}.$
\end{theorem}

\section{Concluding remarks}\label{conc}
While we have stated results for Calabi-Yau hypersurfaces in $\P^n$,
the methods presented here can be extended to handle Calabi-Yau
hypersurfaces in weighted projective spaces, and probably have
applications to hypersurfaces in toric varieties.

One can generalize in another direction by replacing the field of
complex numbers by an arbitrary field of characteristic zero.  There
is an obstacle, however, to using a field of characteristic $p$.  
Namely, our 
construction of the $L_\infty$ structure $\check{Q}$ on the space of
\v Cech cochains requires division by a factorial.
Equation \ref{eq3}, for example, has a factor of $\frac{1}{2}.$  Our
model is defined for large primes, $p>n!$.

As mentioned in the introduction, the moduli space $\M$ is a Frobenius
manifold.  We have described an algebra structure on the tangent
spaces $T_a \M$.  A Frobenius structure implies a great deal more.  
The products on $T_a \M$ are connected
by a \emph{potential} function.  That is, there exist
coordinates on $\M$ so that the
structure constants of the shifted cohomology rings are the third
derivatives of a single function satisfying the WDVV equations. 
We are interested what the special basis of $R$ must be and how 
the potential function appears
in our construction
and will return to this point in a later paper.

%
%

\end{document}